\documentclass[12pt]{amsart}
\usepackage{geometry,amsmath,amssymb,graphicx,amsthm}
\usepackage{cite} 
\usepackage{xspace}
\newcommand{\dueto}[1]{\textup{\textbf{(#1) }}}
\newcommand{\tmem}[1]{{\em #1\/}}
\newcommand{\tmmathbf}[1]{\ensuremath{\boldsymbol{#1}}}
\newcommand{\tmop}[1]{\ensuremath{\operatorname{#1}}}
\newcommand{\tmtextit}[1]{{\itshape{#1}}}
\newcommand{\udots}{{\mathinner{\mskip1mu\raise1pt\vbox{\kern7pt\hbox{.}}\mskip2mu\raise4pt\hbox{.}\mskip2mu\raise7pt\hbox{.}\mskip1mu}}}
\newtheorem{proposition}{Proposition}
\newtheorem{theorem}{Theorem}
\newtheorem{remark}{Remark}

\newtheorem*{definition}{Definition}
\newtheorem*{bc}{Boundary Conditions}

\renewcommand{\tilde}{\widetilde}
\newcommand{\plussign}{$+$\xspace}
\newcommand{\minussign}{$-$\xspace}

\begin{document}

\title{Metaplectic Ice}

\author[B. Brubaker]{Ben Brubaker}
\address{Department of Mathematics \\ MIT \\ Cambridge MA 02139-4307}
\email{brubaker@math.mit.edu}
\author[D. Bump]{Daniel Bump}
\address{Department of Mathematics \\ Stanford University \\ Stanford CA 94305-2125}
\email{bump@math.stanford.edu}
\author[G. Chinta]{Gautam Chinta}
\address{Department of Mathematics \\ The City College of CUNY \\ New York, NY 10031}
\email{chinta@sci.ccny.cuny.edu}
\author[S. Friedberg]{Solomon Friedberg}
\address{Department of Mathematics \\ Boston College \\ Chestnut Hill MA 02467-3806}
\email{friedber@bc.edu}
\author[P. E. Gunnells]{Paul E. Gunnells}
\address{Department of Mathematics and Statistics \\ University of Massachusetts Amherst\\ Amherst, MA 01003}
\email{gunnells@math.umass.edu}

\subjclass[2010]{Primary: 11F68, Secondary: 11F70, 16T25, 22E50}
\keywords{$p$-adic Whittaker function, metaplectic group, Gelfand-Tsetlin
pattern, six-vertex model, statistical mechanics, Yang-Baxter
equation, multiple Dirichlet series}

\begin{abstract}
We study spherical Whittaker functions on a metaplectic cover of
$\tmop{GL}(r+1)$ over a nonarchimedean local field using lattice
models from statistical mechanics. An explicit description of this Whittaker
function was given in terms of Gelfand-Tsetlin patterns in \cite{eisenxtal,
McNamara}, and we translate this description into an expression of the
values of the Whittaker function as partition functions of a six-vertex
model. Properties of the Whittaker function may then be expressed in terms of
the commutativity of row transfer matrices potentially amenable to proof using
the Yang-Baxter equation. We give two examples of this: first, the equivalence
of two different Gelfand-Tsetlin definitions, and second, the effect of
the Weyl group action on the Langlands parameters. The second example is
closely connected with another construction of the metaplectic Whittaker
function by averaging over a Weyl group action \cite{CG, ChintaOffen}.

\end{abstract}

\maketitle

\section{Introduction}

The study of spherical Whittaker functions of reductive groups over
local fields is of fundamental importance in number theory and
representation theory.  Recently, in two separate series of papers,
the authors and their collaborators have studied Whittaker functions
on metaplectic covers of such groups.  The goal of this paper is to
introduce a new method for describing such $p$-adic metaplectic
Whittaker functions: two-dimensional lattice models of statistical mechanics.
In such a model, one defines the partition function to be a weighted sum
over states of the model.  We show that there exists a choice of weights for
which the partition functions are metaplectic Whittaker functions. Baxter
\cite{baxter} developed important techniques for evaluating the partition
functions of lattice models including the so-called ``commutativity of
transfer matrices'' and the use of the Yang-Baxter equation. We discuss how
these methods relate to our descriptions of Whittaker functions and to prior
work.

Two different explicit formulas have been given in
{\cite{ChintaOffen}} and {\cite{McNamara}} for the spherical Whittaker
function on a metaplectic cover of $\tmop{GL} (r+1)$ over a
non-archimedean local field. The first of these is expressed in terms
of a Weyl group action described in {\cite{CG}},
the second in terms of a function on Gelfand-Tsetlin patterns
initially introduced in \cite{wmd3}. In fact, this latter
representation belongs to a family of explicit formulas, one for
each reduced expression of the long element of the Weyl group as a
product of simple reflections.  Two such reduced expressions in type A
are particularly nice, and lead to representations of the Whittaker
function as sums over Gelfand-Tsetlin patterns. In keeping with
earlier works, we refer to the two different descriptions as ``Gamma''
and ``Delta'' rules. The main result of \cite{wmd5book} is a
combinatorial proof that these two definitions are in fact equal. This
equality allows one to prove the analytic properties of an associated
global object (a multiple Dirichlet series) by applying Bochner's
convexity principle.

In the following section we demonstrate that the Gelfand-Tsetlin
patterns we are concerned with are in bijection with admissible states
of the six-vertex model having certain fixed boundary conditions. After
recalling the description of the metaplectic Whittaker function as a function
on Gelfand-Tsetlin patterns in Section 3, we use the bijection with the
lattice model in Section 4 to express both the Gamma and Delta descriptions of
the Whittaker function as partition functions for certain respective choices
of Boltzmann weights.

In Section 5, we take the connection with statistical models
further. We show that the necessary result for demonstrating the
equivalence of the Gamma and Delta descriptions may be reformulated in
terms of the commutativity of transfer matrices.  Baxter \cite{baxter} advocated the
use of the Yang-Baxter equation for demonstrating this commutativity.
In Section 6, we explain how this is carried out in the context of the six-vertex model
and we speculate about the possibility of such an equation in the
metaplectic case.

Finally we discuss the Weyl group action on metaplectic Whittaker
functions, initially established by Kazhdan and Patterson \cite{kp}, which
plays a critical role in the explicit formulas of \cite{ChintaOffen}.  When
the degree of the cover is $1$, i.e. the linear case, the $p$-adic spherical
Whittaker function is essentially a Schur polynomial
by results going back to Shintani \cite{Shintani}.
The Weyl group action is thus closely related to the standard permutation 
action on polynomials in $r+1$ variables.  In \cite{hkice}, this Whittaker
function (or equivalently, the Schur polynomial multiplied by a $q$-deformation
of the Weyl denominator) is realized as a partition 
function on a six-vertex model and its
properties are studied via instances of the Yang-Baxter equation. On
the other hand, as soon as the degree of the cover is greater than
$1$, the action looks rather different
(cf.~\eqref{defn:PQ}--\eqref{eq:fe}).  Nevertheless, we may ask
whether these functional equations may also be phrased in terms of
transfer matrices and a Yang-Baxter equation, and in this final
section, we present evidence towards an affirmative answer.

This work was partially supported by the following grants: NSF grants
DMS-0844185 (Brubaker), DMS-0652817 and DMS-1001079 (Bump), DMS-0847586
(Chinta), DMS-0652609 and DMS-1001326 (Friedberg), DMS-0801214 (Gunnells) 
and NSA grant H98230-10-1-0183 (Friedberg).

\section{\label{bijectionsection}Six Vertex Model and Gelfand-Tsetlin
Patterns}

In this section, we demonstrate a bijection between strict
Gelfand-Tsetlin patterns and admissible states of the six vertex model
(or ``square ice'') on a finite square lattice with certain fixed
boundary conditions. The boundary conditions on ice were known to
Hamel and King, who presented bijections between ice and patterns
related to the symplectic group in \cite{hamelkinguturn}.  A treatment
tailored to the aims of the present paper was given in \cite{hkice},
whose terminology we now recall.

A {\it Gelfand-Tsetlin pattern} of rank $r$ is a triangular array of integers
\begin{equation}
  \mathfrak{T} = \left\{ \begin{array}{ccccccc}
    a_{0, 0} &  & a_{0, 1} & \cdots & a_{0, r - 1} &  & a_{0, r}\\
    & a_{1, 1} &  & \cdots &  & a_{1, r} & \\
    &  & \ddots &  & \udots &  & \\
    &  &  & a_{r, r} &  &  & 
  \end{array} \right\} \label{gtindexing}
\end{equation}
in which the rows interleave: $a_{i - 1, j - 1} \geqslant a_{i, j}
\geqslant a_{i - 1, j}$. The set of all Gelfand-Tsetlin patterns with fixed
top row is in bijection with basis vectors of a corresponding highest weight
representation of $\tmop{GL} (r + 1, \mathbb{C})$. Indeed, any given top row
$(a_{0, 0}, a_{0, 1}, \ldots, a_{0, r})$ is a partition which may be
regarded as a dominant weight of the $\tmop{GL} (r + 1, \mathbb{C})$ weight
lattice. Each successive row of a pattern then records a branching rule down
to a highest weight representation on a subgroup of rank one less. We will
focus mainly on the set of {\it strict} Gelfand-Tsetlin patterns, whose entries in
horizontal rows are strictly decreasing. In terms of representation theory,
these patterns result from branching through strictly dominant highest
weights. Top rows of strict Gelfand-Tsetlin patterns are then indexed by
strictly dominant weights $\lambda + \rho$ where $\lambda$ is a dominant
weight and $\rho$ is the Weyl vector $(r, r - 1, \ldots, 0)$.

Now we come to lattice models. The six-vertex model consists of labelings of edges
in a square grid where each vertex has adjacent edges in one of six admissible
configurations. This model is sometimes referred to as ``square ice'' where
each vertex of the grid represents an oxygen atom and the 6 admissible ways of
labeling adjacent edges correspond to the number of ways in which two of the
four edges include a nearby hydrogen atom. If we represent adjacent hydrogen
atoms by incoming arrows, and locations where there is no adjacent hydrogen
atom by outgoing arrows, the six admissible states are as follows.
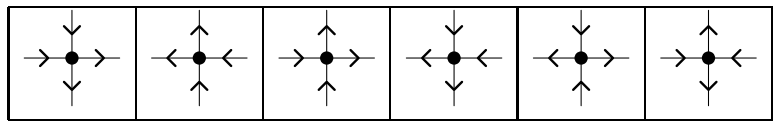
\[ \begin{array}{|c|c|c|c|c|c|c|}
     \hline
     \includegraphics{arrows1.mps} & \includegraphics{arrows2.mps} &
     \includegraphics{arrows3.mps} & \includegraphics{arrows4.mps} &
     \includegraphics{arrows6.mps} & \includegraphics{arrows5.mps}\\
     \hline
   \end{array} \]

We will use a representation consisting of a lattice whose edges are
labeled with signs $+$ or $-$, called {\it spins}. To relate this
to the previous description, interpret a right-pointing or down-pointing
arrow as a $+$, and a left-pointing or up-pointing arrow as $-$. We then
find the following six configurations. (The index $i$ in the table indicates
the row to which the vertex belongs, and will be used in later sections.)
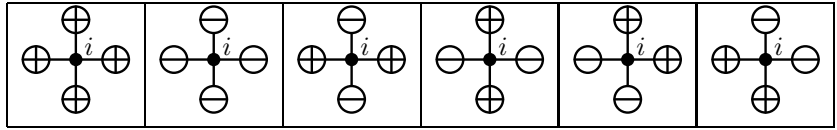
\[ \begin{array}{|c|c|c|c|c|c|c|}
     \hline
     \includegraphics{gamma1a.mps} &
     \includegraphics{gamma6a.mps} & \includegraphics{gamma4a.mps} &
     \includegraphics{gamma5a.mps} & \includegraphics{gamma2a.mps} &
     \includegraphics{gamma3a.mps}\\
     \hline
   \end{array} \]

The rectangular lattices we consider will be finite, with boundary conditions
chosen so that the admissible configurations are in bijection with strict
Gelfand-Tsetlin patterns with fixed rank $r$ and top row $\lambda + \rho$ as
above. Here $\lambda=(\lambda_r,\dots,\lambda_1,\lambda_0)$ with
$\lambda_j \geq \lambda_{j-1}$ for all $j$, and we suppose 
that $\lambda_0=0$.  

\begin{bc} The rectangular grid is to
have $\lambda_r + r + 1$ columns (labeled $0$ through $\lambda_r + r$
increasing from right to left) and $r + 1$ rows. Then with $\lambda + \rho =
(\lambda_r + r, \lambda_{r - 1} + r - 1, \ldots, 0)$, we place a \minussign spin at
the top of each column whose label is one of the distinct parts of $\lambda +
\rho$, i.e. at columns labeled $\lambda_j + j$ for $0 \leqslant j \leqslant
r$. We place a \plussign spin at the top of each of the remaining columns.
Furthermore, we place a \plussign spin at the bottom of every column and on
the left-hand side of each row and a \minussign spin on the right-hand side of
each row.
\end{bc}

For example, put $r =
2$, and take $\lambda = (3, 2, 0)$, so that $\lambda + \rho = (5, 3,0)$.  
Then we have the following boundary conditions for the ice:
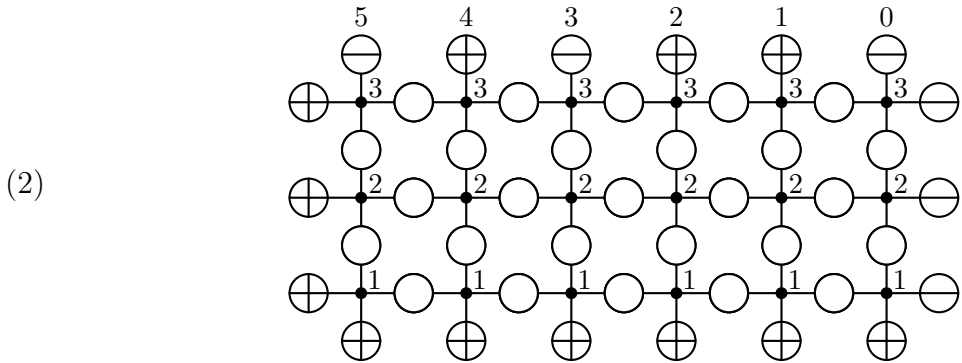
\begin{equation} \vcenter{\hbox to 3in{\includegraphics[scale=0.9]{gamma_ice2.mps}}} \label{bdconds} 
\end{equation}
The column labels are written above each column, and row labels have been
placed next to each vertex. These row labels will be used in Section~\ref{sec:imwf},
but need not concern us now. 
The edge spins have been placed inside circles
located along the boundary. The remaining open circles indicate interior spins
not determined by our boundary conditions, though any filling of the grid must
use only the six admissible configurations in the above table. Such an admissible 
filling of the finite lattice having above boundary conditions will be referred to as a
\tmtextit{state} of ice.

\begin{proposition}
  \label{bijprop}Given a fixed rank $r$ and a dominant weight $\lambda
  = (\lambda_r, \ldots, \lambda_1, 0)$, there is a bijection between strict
  Gelfand-Tsetlin patterns with top row $\lambda + \rho$ and admissible states
  of ice having boundary conditions determined by $\lambda$ as above.
\end{proposition}

\begin{proof}
  We begin with a strict Gelfand-Tsetlin pattern.
  Each row of the Gelfand-Tsetlin pattern will correspond to the set of spins
  located between numbered rows of ice, the so-called ``vertical spins''
  since they lie on vertical edges of the grid, as follows. To each entry $a_{i, j}$ in
  the Gelfand-Tsetlin pattern, we assign a \minussign to the vertical spin between rows labeled $r+2-i$
  and $r+1-i$ in the column labeled $a_{i, j}$. (Recall that we are using
  decreasing row labels from top to bottom as in the example (\ref{bdconds}).)
  The remaining vertical spins are assigned \plussign.
  
  It remains to assign horizontal spins, but these are already uniquely
  determined since the left and right edge horizontal spins have been assigned
  and each admissible vertex configuration has an even number of adjacent \plussign
  spins. We must only verify that the resulting configuration uses only the 6
  admissible configurations (from the 8 having an even number of \plussign signs) for
  the corresponding ice. This is easily implied by the interleaving condition on entries 
  in the Gelfand-Tsetlin pattern, which is violated if one of the two inadmissible configurations
  appears. See Lemma 2 of {\cite{hkice}} for more details.
\end{proof}

A simple example illustrates the bijection:
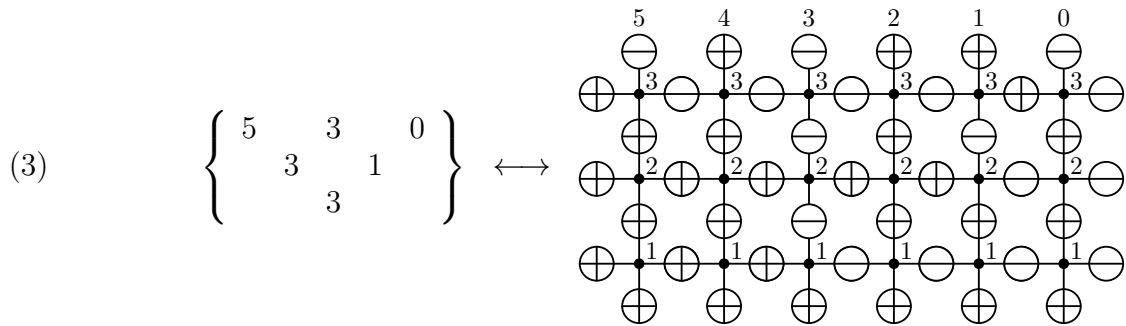
\begin{equation}
  \begin{array}{ccc}
    \mathfrak{} \left\{ \begin{array}{ccccc}
      5 &  & 3 &  & 0\\
      & 3 &  & 1 & \\
      &  & 3 &  & 
    \end{array} \right\} & \longleftrightarrow &
    \vcenter{\hbox to 2in{\includegraphics[scale=0.8]{gamma_ice1.mps}}}
  \end{array} \label{workingexample}
\end{equation}

\section{\label{whitdefined}Metaplectic Whittaker Functions and Patterns}

We now discuss the relation between the spherical metaplectic
Whittaker function on the $n$-fold cover of of $\tmop{GL} (r + 1)$
over a non-archimedean local field and Gelfand-Tsetlin patterns.
Such a relationship, described globally, was conjectured in
{\cite{wmd3}} and was established in {\cite{eisenxtal}}.  Though it is
possible to pass from the global result to its local analogue
(cf. \cite{Friedberg-McNamara}), a direct local proof was given by
McNamara \cite{McNamara}, expressing a metaplectic spherical Whittaker
function as a generating function supported on strict Gelfand-Tsetlin
patterns. In this section, we recall two formulations of this explicit
description, following \cite{wmd5book}. In Section~\ref{sec:imwf}, we
will explain their translations to square ice 
via the bijection of Proposition \ref{bijprop}.

Let $\tilde{G(F)}$ denote the $n$-fold metaplectic cover of $G(F) =
\tmop{GL}(r+1,F)$, where $F$ denotes a nonarchimedean local field having ring of integers
$\mathfrak{o}_F$ and residue field of order $q$. (There are several
related such extensions, but we will use the one in~\cite{kp} where
their parameter $c=0$.) We assume $2n$ divides $q-1$, which guarantees that
$F$ contains the group $\mu_{2n}$ of $2n$th roots of unity. The group
$\widetilde{G(F)}$ is a central extension of $G(F)$ by $\mu_n$:
$$ 1 \longrightarrow \mu_n \rightarrow \widetilde{G(F)}
\stackrel{\pi}{\longrightarrow} G(F) \longrightarrow 1. $$
We will identify $\mu_n\subset F$ with the group $\mu_n\subset\mathbb{C}$
of complex $n$-th roots of unity by some fixed isomorphism. For
convenience, we will sometimes denote $\widetilde{G(F)}$ as just
$\widetilde{G}$, and if $H$ is an algebraic subgroup of $G$, we
may denote by $\widetilde{H}$ the preimage of $H(F)$ in $\widetilde{G}$.

For details of the construction of the metaplectic group and results about its
representations, see~\cite{mcnamarafoundations} in this volume. Let $B(F)$
be the standard Borel subgroup of upper triangular matrices in $G(F)$,
and let $T(F)$ be the diagonal maximal torus. Then $B(F)=T(F) U(F)$ where
$U(F)$ is the unipotent radical of $B(F)$. The metaplectic cover splits
over various subgroups of $G(F)$; for us it is relevant that it splits
over $U(F)$ and over $K := G(\mathfrak{o}_F)$, the standard maximal compact
subgroup. By abuse of notation, we will sometimes denote by $K$ the
homomorphic image of $K$ in $\widetilde{G}$ under this splitting.

Let $\mathbf{s}:G(F)\to\widetilde{G}$ be any map such that $\pi\circ\mathbf{s}$
is the identity map on $G(F)$. Then the map $\sigma:G(F)\times G(F)\to\mu_n$ such
that 
\[\mathbf{s}(g_1)\mathbf{s}(g_2)=\sigma(g_1,g_2)\mathbf{s}(g_1g_2)\]
is a 2-cocycle defining a class in $H^2(G(F),\mu_n)$. A particular
such cocycle was considered by Matsumoto~\cite{Matsumoto}, Kazhdan and
Patterson~\cite{kp} and Banks, Levy and Sepanski~\cite{BanksLevySepanski}. 
By these references, such as~\cite{kp} Section~0.1, the
map $\mathbf{s}$ may be chosen so that the restriction of $\sigma$ to $T(F)$ is
given by the formula
\begin{equation}
\label{cocycle}
  \sigma\left(\left(\begin{array}{ccc}t_1\\&\ddots\\&&t_{r+1}\end{array}\right),
  \left(\begin{array}{ccc}u_1\\&\ddots\\&&u_{r+1}\end{array}\right)\right)
  =\prod_{i<j}(t_i,u_j)_n.
\end{equation}
The cocycle $\sigma$ also has the property that $\sigma(u,g)=\sigma(g,u)=1$
if $u\in U(F)$, and so the restriction of $\mathbf{s}$ to $U(F)$ is a
homomorphism to $\widetilde{G}$.

We will call a representation $\pi$ of $\tilde{G}$ or any subgroup {\it genuine}
if $\pi(\zeta g)=\zeta\pi(g)$ when $\zeta\in\mu_n$. Recall that
$\widetilde{T}$ denotes the inverse image under $\pi$ of the maximal torus
$T(F).$ The center $Z(\widetilde{T})$ of $\widetilde{T}$ has finite
index. Since $\widetilde{T}$ is a Heisenberg group, characters of
$Z(\widetilde{T})$ parametrize the irreducible genuine representations of
$\widetilde{T}$ as follows. Let $\chi$ be a quasi-character of
$Z(\widetilde{T})$ that is unramified. This means that it is trivial on the
maximal compact subgroup of this abelian group. Let $A$ be a maximal abelian
subgroup of $\widetilde{T}$. Extend $\chi$ to $A$, then induce it to
$\widetilde{T}$. By Clifford theory, the resulting representation $i(\chi)$ of
$\widetilde{T}$ is genuine, irreducible, and finite-dimensional. It does not
depend on either $A$ or the extension of $\chi$ to $A$.

We extend $i(\chi)$ to the inverse image $\widetilde{B}$ of $B(F)$ in such a
way that $\mathbf{s}(U(F))$ acts trivially. We then consider the
representation of $\widetilde{G}$ obtained by normalized induction. We call
the vector space of the resulting representation $I(\chi)$. It has a
one-dimensional space of $K$-fixed, i.e.~spherical, vectors.

Let $\psi_0: F \rightarrow \mathbb{C}$ be an additive character
that is trivial on $\mathfrak{o}_F$ but on no larger fractional ideal.  Given
a spherical vector $f$ in $I(\chi)$, we may construct the spherical
metaplectic Whittaker function via the integral
\begin{equation} 
W(g) = \int_{U(F)} f(w_0 \mathbf{s}(u) g) \psi(u) du, \label{whittakerfunction} 
\end{equation}
where $\psi$ is the character of $U(F)$ given by
$$ 
\psi \left( \begin{pmatrix} 1 & x_{1,2} & \cdots & & x_{1,n} \\ & 1 &
x_{2,3} & \cdots & x_{2,n} \\ & & \ddots & & \vdots \\ & & & & 1
\end{pmatrix} \right) = \psi_0(\sum_i x_{i,i+1})
$$
and $w_0$ is a representative in $K$ of the long Weyl group
element. Strictly speaking, $W(g)$ as we have defined it is an
$i(\chi)$-valued function and should be composed with a natural
choice of linear functional on $i(\chi)$ to obtain a complex-valued
function. For brevity we will sometimes refer to $W(g)$ as the {\it
metaplectic Whittaker function}.

Note that $\chi$ on $Z(\widetilde{T})$ may be described by a choice of $r$
complex numbers $\tmmathbf{s} = (s_1, \ldots, s_r)$. The
transformation property $W(\mathbf{s}(u)gk) = \psi(u) W(g)$ for all $u \in U(F), k \in K$,
implies that it suffices to determine $W$ on the inverse image of
the torus $T(F)$. Moreover, since $W$ is genuine, it is sufficient to specify $W$ on
$\mathbf{s}(T(F))$. Given $\lambda = \sum_i \lambda_i \omega_i$, where
$\omega_i$ are fundamental weights, let $t_\lambda$ be the element of the
torus $T(F)$ given by
$$ t_\lambda = \begin{pmatrix} p^{\lambda_1+\lambda_2+\cdots+\lambda_r} & & &
  &  \\ & p^{\lambda_2+\cdots+\lambda_r}  & & & \\ & & \ddots & \\  & & &
  p^{\lambda_r} & \\ & & & & 1 \end{pmatrix} \quad \text{where $p$ is a
  uniformizer for $\mathfrak{o}_F$.} $$
Let $\mathbf{t}_\lambda=\mathbf{s}(t_\lambda)$.
Due to our assumption that $F$ contains the $2n$-th roots of unity,
$(p,p)_n=1$ and by (\ref{cocycle}) it follows that
$\mathbf{t}_{\lambda+\mu}=\mathbf{t}_\lambda\mathbf{t}_\mu$.

It is not hard to show that $W(\mathbf{t}_\lambda) = 0$ unless $\lambda$ is a dominant weight.

Given any dominant weight $\lambda$, the metaplectic Whittaker function
$W(\mathbf{t}_\lambda)$ may thus be expressed as a Dirichlet series in $r$
complex variables $\tmmathbf{s} = (s_1, \ldots, s_r)$ of the form
\begin{equation}
  \label{eq:aa}
 Z ( \tmmathbf{s} ; \lambda) = \sum_{\tmmathbf{k} = (k_1, \ldots, k_r) \in
   \mathbb{N}^r} H (p^{k_1}, \ldots, p^{k_r} ; \lambda) \, q^{k_1 (1 - 2s_1) +
   \cdots + k_r (1-2 s_r)},  
\end{equation}
where $q$
denotes the cardinality of the residue field $\mathfrak{o}_F / \mathfrak{p}$.
We now describe the function $H (p^{k_1}, \ldots, p^{k_r} ; \lambda)$. See
{\cite{McNamara}} for further details.

The positive integer $n$ will continue to denote the degree of the
metaplectic cover. We define the Gauss sum
\[ g (a, b) = \int_{\mathfrak{o}_F^\times} \left( u, p \right)_n^b \psi_0 \left( p^{a-b} u \right) \, du, \]
where $\left( \cdot, \cdot \right)_n$ denotes the $n$-th power Hilbert
symbol, and we normalize the Haar measure so that $\mu(\mathfrak{o}_F) =
1$. As a further shorthand, for any positive integer $b$, we set
\begin{equation}
  g (b) = g (b - 1, b), \hspace{2em} h (b) = g (b, b) .
  \label{gandhforshort}
\end{equation}
Note that for a fixed base field $F$, these values depend only on $b$
mod $n$. If $n$ divides $b$, in particular if $n=1$, we have
\begin{equation} g(b) = -\frac{1}{q}, \quad h(b) = 1 -
  \frac{1}{q}. \label{nequalone} \end{equation}
We caution the reader that the $q$-powers that appear in the $g(a,b)$
are normalized differently than in the previous works
\cite{eisenxtal, wmd5book, wmd3}; these are the functions denoted $g^\flat$
and $h^\flat$ in \cite{wmd5book}. The function $h$ is a degenerate Gauss sum
whose values may be made explicit, while (if $n\nmid b$) $g(b)$ is a ``true''
Gauss sum.

Any strict Gelfand-Tsetlin pattern $\mathfrak{T}$ with entries indexed as in
(\ref{gtindexing}), we associate a weighting function $\gamma$ to each entry
$a_{i, j}$ with $i \geq 1$ as follows:
\begin{equation}
  \gamma (a_{i, j}) = \left\{ \begin{array}{ll}
    g (b_{i, j}) & \tmop{if} a_{i, j} = a_{i - 1, j - 1},\\
    h (b_{i, j}) & \tmop{if} a_{i - 1, j} \neq a_{i, j} \neq a_{i - 1, j -
    1},\\
    1 & \tmop{if} a_{i, j} = a_{i - 1, j},
  \end{array} \right. \hspace{1em} \tmop{where} \hspace{1em} b_{i, j} =
  \sum^r_{l = j} (a_{i, l} - a_{i - 1, l}). \label{gammagtwt}
\end{equation}
Then we define
\begin{equation}
  G^{\Gamma} ( \mathfrak{T}) = \prod_{i = 1}^r \prod^r_{j = i} \gamma (a_{i,
  j}). \label{Ggamma}
\end{equation}
If $\mathfrak{T}$ is a Gelfand-Tsetlin pattern that is not strict, we define
$G^{\Gamma} ( \mathfrak{T})=0$. We also define
\begin{equation}
  \tmmathbf{k}^{\Gamma} ( \mathfrak{T}) = (k_1^{\Gamma} ( \mathfrak{T}), \ldots,
  k_r^{\Gamma} ( \mathfrak{T})) \hspace{1em} \tmop{where} \hspace{1em}
  k_i^{\Gamma} ( \mathfrak{T}) = \sum^r_{l = i} a_{i, l} - a_{0, l} .
  \label{kgamma}
\end{equation}
In particular, note that both $G^{\Gamma}$ and $\tmmathbf{k}^{\Gamma}$ are defined using
differences of elements above and to the right of $a_{i, j}$. The superscript
$\Gamma$ may be regarded as indicator that these quantities are defined using such ``right-hand''
differences. 

We present these definitions in this {\it ad-hoc} fashion in order to give a brief and 
self-contained treatment, but in fact they have very natural descriptions when reinterpreted as
functions on a Kashiwara crystal graph. See {\cite{wmd5book}} for an extensive
discussion.

As an example, consider the Gelfand-Tsetlin pattern $\mathfrak{T}$ in
(\ref{workingexample}).  Then
\begin{equation}
  (b_{1, 1}, b_{1, 2}, b_{2, 2}) = (1, 1, 2) \hspace{1em} \tmop{so}
  \tmop{that} \hspace{1em} G^{\Gamma} ( \mathfrak{T}) = h (1) g(2),
  \hspace{1em} \tmop{and} \hspace{1em} (k_1, k_2) = (1, 3) .
  \label{gofexample}
\end{equation}

\begin{theorem}\label{thm:bbfm}
  {\dueto{Brubaker, Bump, and Friedberg~{\cite{eisenxtal}};
  McNamara~{\cite{McNamara}}}}Given a dominant weight $\lambda$ and a fixed
  $r$-tuple of non-negative integers $\tmmathbf{k} = (k_1, \ldots, k_r)$, the
  function \ $H (p^{k_1}, \ldots, p^{k_r} ; \lambda)$ appearing in the p-adic
  Whittaker function $W(\mathbf{t}_\lambda)$ is given by
  \[ H (p^{\tmmathbf{k}} ; \lambda) : = H (p^{k_1}, \ldots, p^{k_r} ; \lambda)
     = \sum_{\tmmathbf{k}^{\Gamma} ( \mathfrak{T}) = \tmmathbf{k}} G^{\Gamma} (
     \mathfrak{T}_{}), \]
  where the sum is over all Gelfand-Tsetlin patterns with top row
  corresponding to $\lambda + \rho$ satisfying the subscripted condition.
\end{theorem}

There is a second explicit description of $H (p^{\tmmathbf{k}} ; \lambda)$ in
terms of ``left-hand'' differences using  functions $G^{\Delta}$ and
$\tmmathbf{k}^{\Delta}$ that are analogous to those defined in (\ref{Ggamma}) and (\ref{kgamma})
respectively. Assuming that $\mathfrak{T}$ is strict, set
\begin{equation}
  \delta (a_{i, j}) = \left\{ \begin{array}{ll}
    g (c_{i, j}) & \tmop{if} a_{i, j} = a_{i - 1, j},\\
    h (c_{i, j}) & \tmop{if} a_{i - 1, j} \neq a_{i, j} \neq a_{i - 1, j -
    1},\\
    1 & \tmop{if} a_{i, j} = a_{i - 1, j - 1},
  \end{array} \right. \hspace{1em} \tmop{where} \hspace{1em} c_{i, j} =
  \sum^j_{l = 1} (a_{i - 1, l - 1} - a_{i, l}) \label{deltagtwt}
\end{equation}
and define
\begin{equation}
  G^{\Delta} ( \mathfrak{T}) = \prod_{i = 1}^r \prod^r_{j = i} \delta (a_{i,
  j}) . \label{Gdelta}
\end{equation}
If $\mathfrak{T}$ is not strict, define $G^\Delta(\mathfrak{T})=0$. We also set
\begin{equation}
  \tmmathbf{k}^{\Delta} ( \mathfrak{T}) = (k_1^{\Delta} ( \mathfrak{T}), \ldots,
  k_r^{\Delta} ( \mathfrak{T})) \hspace{1em} \tmop{where} \hspace{1em}
  k_i^{\Delta} ( \mathfrak{T}) = \sum^i_{l = 1} a_{0, l-1} - a_{r + 1 - i, r + 1
  - l} . \label{kdelta}
\end{equation}
The main theorem of {\cite{wmd5book}} is as follows:

\begin{theorem}
  \label{statementA}{\dueto{Statement A of Brubaker, Bump and Friedberg {\cite{wmd5book}}}}
  Given a dominant weight $\lambda$ and a fixed
  $r$-tuple of non-negative integers $\tmmathbf{k} = (k_1, \ldots, k_r)$,
  \begin{equation}
    \sum_{k^{\Gamma} ( \mathfrak{T}) = \tmmathbf{k}} G^{\Gamma} (
    \mathfrak{T}_{}) = \sum_{k^{\Delta} ( \mathfrak{T}) = \tmmathbf{k}}
    G^{\Delta} ( \mathfrak{T}_{}), \label{gammaisdelta}
  \end{equation}
where the sums each run over all Gelfand-Tsetlin patterns with top row
  corresponding to $\lambda + \rho$ satisfying the subscripted condition.
\end{theorem}

As an immediate corollary, we have a second description of the $p$-adic Whittaker
function in terms of $G^{\Delta}$ and $\tmmathbf{k}^{\Delta}$. We refer to these two
recipes on the left- and right-hand sides of (\ref{gammaisdelta}) as the
$\Gamma$- and $\Delta$-rules, respectively.  

In fact, there are many other descriptions for the Whittaker function,
though these are generally much more difficult to write down as explicitly.  Indeed, as
explained in {\cite{BerensteinZelevinsky, Littelmann}},
there exist bases for highest weight representations
corresponding to any reduced expression for the long element $w_{0}$
of the Weyl group of $\tmop{GL} (r + 1)$ --- $S_r$, the symmetric
group on $r$ letters --- as a product of simple reflections
$\sigma_i$.  These make use of the Kashiwara crystal graph and are commonly
called {\it string bases}. Using these bases, one may make a correspondence between long
words and recipes for the Whittaker function (cf. {\cite[Chapter
2]{wmd5book}}).  From this perspective, the $\Gamma$-rule corresponds
to the word
\[ w_0 = \sigma_1 (\sigma_2 \sigma_1) \cdots (\sigma_r \sigma_{r - 1} \cdots
   \sigma_1), \]
whereas the $\Delta$-rule corresponds to the word
\[ w_0 = \sigma_r (\sigma_{r - 1} \sigma_r) \cdots (\sigma_1 \sigma_2 \cdots
   \sigma_r) . \]
These two words are as far apart as possible in the lexicographic ordering of
all reduced decompositions. The proof of Theorem \ref{statementA} as given in 
\cite{wmd5book} uses a blend of combinatorial
arguments to give various equivalent forms of the identity (\ref{gammaisdelta}) as
we move through the space of long words. 
We highlight various aspects of the proof in more detail now.

The proof is by induction on the rank $r$. The inductive hypothesis
allows us to equate any two recipes for the Whittaker function whose
associated long words differ by a sequence of relations obtained from
a lower rank case. For example, assuming the rank 2 case allows us to
perform a braid relation $\sigma_1 \sigma_2 \sigma_1 = \sigma_2
\sigma_1 \sigma_2$, which could be applied to the word corresponding
to the $\Gamma$-rule above. After a series of such identities, we
arrive at two descriptions for the Whittaker function as a weighted
sum over Gelfand-Tsetlin patterns that agree on the bottom $r - 2$
rows of the pattern. Thus, we may restrict our attention to the top
three rows of a rank $r$ pattern. We refer to such three-row arrays of
interleaving integers, where we fix both the top and bottom of the
three rows, as ``short Gelfand-Tsetlin patterns'' and re-index the
three rows as follows:
\begin{equation}
  \mathfrak{t =} \left\{ \begin{array}{ccccccccc}
    \ell_0 &  & \ell_1 &  & \cdots &  & \ell_{r - 1} &  & \ell_r\\
    & a_1 &  & a_2 &  & a_{r - 1} &  & a_r & \\
    &  & m_1 &  & m_2 &  & m_{r - 1} &  & 
  \end{array} \right\} \label{shortgtindexing} .
\end{equation}
These two recipes for the Whittaker function will be called $G^{\Gamma
\Delta}$ (as this recipe uses a right-hand rule for the entries $a_i$ and a left-hand
rule for the entries $m_j$) and $G^{\Delta \Gamma}$ (where the use of rules is
reversed). To be exact, using the definitions in (\ref{gammagtwt}) and
(\ref{deltagtwt}), we have
\[ G^{\Gamma \Delta} ( \mathfrak{t}) = \prod_{i = 1}^r \gamma (a_i) \prod^{r -
   1}_{j = 1} \delta (m_j), \hspace{1em} \tmop{and} \hspace{1em} G^{\Delta
   \Gamma} ( \mathfrak{t}) = \prod_{i = 1}^r \delta (a_i) \prod^{r - 1}_{j =
   1} \gamma (m_j) . \]
Rather than define functions $k^{\Gamma \Delta}$ and $k^{\Delta \Gamma}$ on
short patterns in analogy to the recipes above, it is enough to specify the
middle row sum as the other rows are fixed.

Before stating the reduction, we require one final ingredient. There is a
natural involution $q_r$ on short Gelfand-Tsetlin patterns of rank $r$, given
by acting on middle row entries according to
\[ q_r : \hspace{0.75em} a_i \longmapsto \max (\ell_{i - 1}, m_{i - 1}) + \min
   (\ell_i, m_i) - a_i = : a_i', \]
where if $i = 0$ we understand that $\max (\ell_0, m_0) = \ell_0$ and if $i =
r$, $\min (\ell_r, m_r) = \ell_r$. This involution $q_r$ is used by Berenstein
and Kirillov (cf.~\cite{BerensteinKirillov}) to define a Sch\"utzenberger
involution on Gelfand-Tsetlin patterns. Brubaker, Bump and Friedberg use the involution $q_r$ to give
the following reduction of Statement A in {\cite{wmd5book}}.

\begin{theorem}
  \label{statementB}{\dueto{Brubaker, Bump and Friedberg; Statement B of {\cite{wmd5book}}}} Fix an $(r +
  1)$-tuple of positive integers $\tmmathbf{\ell} = (\ell_0, \ldots, \ell_r)$,
  an $(r - 1)$-tuple of positive integers $\tmmathbf{m} = (m_1, \ldots, m_{r -
  1}),$ and a positive integer $k.$ Then
  \[ \sum_{\sum a_i = k} G^{\Gamma \Delta} ( \mathfrak{t}) = \sum_{\sum a_i' =
     k'} G^{\Delta \Gamma} (q_r (\mathfrak{t})), \]
  where $a_i'$ are the entries of $q_r ( \mathfrak{t})$, $k' = \sum_i \ell_i +
  \sum_j m_j - k$, and the sums range over all short patterns with top row
  $\tmmathbf{\ell}$ and bottom row $\tmmathbf{m}$ satisfying the indicated
  condition.
\end{theorem}

See {\cite{gelbvol}} and Chapter 6 of {\cite{wmd5book}} for a full proof of the
reduction from Statement A to Statement B. As noted above, the proof of
Statement B proceeds through a series of additional reductions which occupy 
thirteen chapters of {\cite{wmd5book}}. In brief, for ``generic'' short patterns $\mathfrak{t}$,
the Sch\"utzenberger involution $q_r$ gives a finer equality $G^{\Gamma
\Delta} ( \mathfrak{t}) = G^{\Delta \Gamma} (q_r (\mathfrak{t}))$, which implies
the equality of Statement B summand by summand.  
By ``generic'' we mean that the entries of the short pattern are in general
position -- in particular, for all $i$, $\ell_i \neq m_i$ using the notation 
of~(\ref{shortgtindexing}). Note that the
Sch\"utzenberger involution does not necessarily preserve strictness for all
patterns in the remaining non-generic cases, and one needs much more subtle
arguments to handle these short patterns.  For such patterns, Statement B
is not in fact true summand by summand, and one does need
to sum over all short patterns with fixed row sum to obtain equality.

As an alternative to establishing Statement B, we mention
that one could also prove Theorem \ref{statementA}
by computing the Whittaker integral in two ways, mimicking the techniques of
{\cite{eisenxtal}}, thus obtaining a proof via decomposition theorems in
algebraic groups which respect the metaplectic cover. In subsequent sections
of this paper, 
we propose a third way of viewing these theorems using ice-type models, which
portends new connections between number theory/representation theory and
statistical physics.

\section{\label{sec:imwf}Ice and Metaplectic Whittaker Functions} \ \ \ \ \ \

In statistical mechanics, one attempts to infer global behavior from local
interactions. In the context of lattice models, this means that we attach a
\tmtextit{Boltzmann weight} to each vertex in the grid and for each
admissible state of the model, we consider the product of all Boltzmann
weights ranging over all vertices of the grid. Then one can attempt to
determine the \tmtextit{partition function} of the lattice model, which is
simply the sum over all admissible states of the associated weights. In this
section, we explain how to obtain the metaplectic spherical Whittaker function
of Section \ref{whitdefined} as the partition function of a lattice model with
boundary conditions as defined in Section \ref{bijectionsection}.

We make use of the two sets of Boltzmann weights $B^\Gamma$ and $B^\Delta$. When
these weights are applied to an admissible state of ice, we refer to the resulting configuration
as {\it Gamma ice} or {\it Delta ice}, respectively. In order to indicate which set of weights is being
used at a particular vertex, we use $\bullet$ for Gamma ice and $\circ$ for Delta ice.
\begin{equation}
  \begin{array}{|l|l|l|l|l|l|l|}
    \hline
    \begin{array}{c}
      \tmop{Gamma}\\
      \text{Ice}
    \end{array} & \vcenter{\vspace{3pt}\hbox to 40pt{\includegraphics{gamma1a.mps}}\vspace{3pt}} &
    \vcenter{\hbox to 40pt{\includegraphics{gamma6a.mps}}} & \vcenter{\hbox
      to 40pt{\includegraphics{gamma4a.mps}}} &
    \vcenter{\hbox to 40pt{\includegraphics{gamma5a.mps}}} & \vcenter{\hbox
      to 40pt{\includegraphics{gamma2a.mps}}} &
    \vcenter{\hbox to 40pt{\includegraphics{gamma3a.mps}}}\\
    \hline
    \text{\begin{tabular}{c}
      Boltzmann\\
      weight $B^{\Gamma}$ 
    \end{tabular}} & 1 & z_i & g (a) & z_i & h (a) z_i & 1\\
    \hline
    \begin{array}{c}
      \tmop{Delta}\\
      \text{Ice}
    \end{array} & \vcenter{\vspace{3pt}{\hbox to 40pt{\includegraphics{delta5a.mps}}}\vspace{3pt}} &
    \vcenter{\hbox to 40pt{\includegraphics{delta4a.mps}}} & \vcenter{\hbox
      to 40pt{\includegraphics{delta6a.mps}}} &
    \vcenter{\hbox to 40pt{\includegraphics{delta1a.mps}}} & \vcenter{\hbox
      to 40pt{\includegraphics{delta2a.mps}}} &
    \vcenter{\hbox to 40pt{\includegraphics{delta3a.mps}}}\\
    \hline
    \text{\begin{tabular}{c}
      Boltzmann\\
      weight $B^{\Delta}$
    \end{tabular}} & 1 &  g(a) z_i  & 1 & z_i & h(a) z_i & 1 \\
    \hline
  \end{array}
\label{weighttable}
\end{equation}

In giving these Boltzmann weights, we have made use of the notation in
(\ref{gandhforshort}).  For Gamma ice, the constant $a$
equals the number of \plussign signs in the $i$-th row to the right of
the vertex $\bullet$. For Delta ice, the constant $a$ equals the number
of \minussign signs in the $i$-th row to the left of the vertex
$\circ$. In either case, we refer to this constant as the ``charge''
at the vertex. Note by our definitions in (\ref{gandhforshort}), the
Boltzmann weights only depend on the charge mod $n$.  The weights
$B^{\Gamma}$ and $B^{\Delta}$ also depend on parameters $z_i$, where
$i$ indicates the row in which the vertex is found.  
For Gamma ice, the row numbers decrease from $r+1$ to $1$
as we move from top to bottom as in the example (\ref{bdconds}), while for Delta ice, the row numbers increase from $1$ to $r+1$.
These $z_i$ are
referred to as ``spectral parameters.'' We often suppress the dependence
of $B^{\Gamma}$ and $B^{\Delta}$ on the spectral parameters $z_i$,
$1\leq i\leq r+1$.  Let $\tmmathbf{z}=(z_1,\dots,z_{r+1})$.

Given an admissible state of Gamma ice (or Delta ice, respectively) $\mathcal{S}$, we define
\begin{equation}
  \mathcal{G}^{\Gamma} ( \mathcal{S},\tmmathbf{z}) = \prod_{v \in \mathcal{S}} B^{\Gamma}
  (v), \hspace{1em} \mathcal{G}^{\Delta} ( \mathcal{S},\tmmathbf{z}) = \prod_{v \in
  \mathcal{S}} B^{\Delta} (v), \label{icygs}
\end{equation}
where the product (in either case) is taken over all vertices in the state of ice $\mathcal{S}$.

\begin{proposition}\label{matching}
  Under the bijection of Proposition \ref{bijprop}, with strict pattern
  $\mathfrak{T}$ corresponding to an admissible state of Gamma ice
  $\mathcal{S_{}}$, then $G^{\Gamma} ( \mathfrak{T})$ as defined in
  (\ref{Ggamma}) is related to $\mathcal{G}^{\Gamma} ( \mathcal{S},\tmmathbf{z})$ in
  (\ref{icygs}) as follows:
  \[ \mathcal{G}^{\Gamma} ( \mathcal{S},\tmmathbf{z}) = G^{\Gamma} ( \mathfrak{T}) z_{r+1}^{d_0
     ( \mathfrak{T}) - d_1 ( \mathfrak{T})} z_r^{d_1 ( \mathfrak{T}) - d_2 (
     \mathfrak{T})} \cdots z_2^{d_{r - 1} ( \mathfrak{T)} - d_r (
     \mathfrak{T})} z_{1}^{d_r ( \mathfrak{T})}, \]
  where $d_i ( \mathfrak{T})$ is the sum of the entries in the $i$-th row of
  the pattern $\mathfrak{T}$.
  
  Similarly, for an admissible state of Delta ice $\mathcal{S}$, $G^{\Delta} (
  \mathfrak{T})$ as defined in (\ref{Gdelta}) is related to
  $\mathcal{G}^{\Delta} ( \mathcal{S},\tmmathbf{z})$ in (\ref{icygs}) by
  \[ \mathcal{G}^{\Delta} ( \mathcal{S}, \tmmathbf{z}) = G^{\Delta} ( \mathfrak{T}) z_{1}^{d_0
     ( \mathfrak{T}) - d_1 ( \mathfrak{T})} z_2^{d_1 ( \mathfrak{T}) - d_2 (
     \mathfrak{T})} \cdots z_r^{d_{r - 1} ( \mathfrak{T)} - d_r (
     \mathfrak{T})} z_{r+1}^{d_r ( \mathfrak{T})} . \]
\end{proposition}

We first illustrate this for Gamma ice with our working example from
(\ref{workingexample}) in rank~2. The admissible Gamma ice $\mathcal{S}$ and
its associated Boltzmann weights are pictured below.

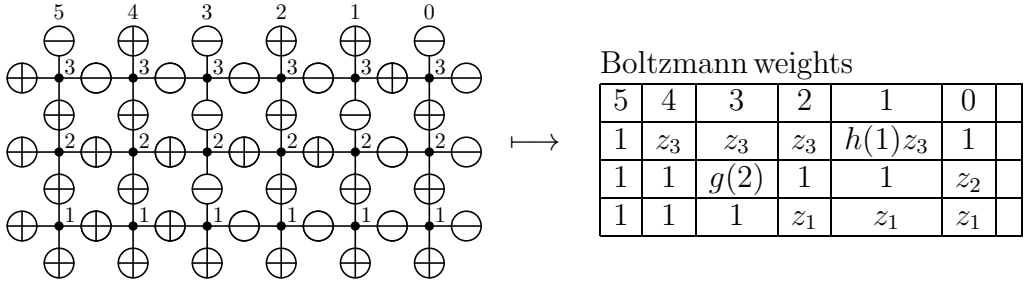
\[ \begin{array}{lll}
     \vcenter{\hbox to 2.5in{\includegraphics[width=2.5in]{gamma_ice1.mps}}} & \longmapsto & \begin{array}{l}
       \tmop{Boltzmann} \tmop{weights}\\
       \begin{array}{|c|c|c|c|c|c|c|}
         \hline
         5 & 4 & 3 & 2 & 1 & 0 & \\
         \hline
         1 & z_3 & z_3 & z_3 & h (1) z_3 & 1 & \\
         \hline
         1 & 1 & g (2) & 1 & 1 & z_2 & \\
         \hline
         1 & 1 & 1 & z_1 & z_1 & z_1 & \\
         \hline
       \end{array}
     \end{array}
   \end{array} \]
Taking the product over all these weights, we obtain $\mathcal{G}^{\Gamma} (
\mathcal{S},\tmmathbf{z}) = h (1) g (2) z_3^4 z_2 z_1^3$, which indeed matches
$G^{\Gamma} ( \mathfrak{T}) z_3^{d_0 ( \mathfrak{T}) - d_1 ( \mathfrak{T})}
z_2^{d_1 ( \mathfrak{T}) - d_2 ( \mathfrak{T})} z_1^{d_2 ( \mathfrak{T})}$
with $\mathfrak{T}$ as in (\ref{workingexample}) and $G^{\Gamma} (
\mathfrak{T})$ as in (\ref{gofexample}).

\begin{remark} \label{remark} The relevant terms in the metaplectic Whittaker function take the form
\begin{equation} G^\Gamma(\mathfrak{T}) q^{k_1^\Gamma(\mathfrak{T}) (1-2s_1) + \cdots + k_r^\Gamma (\mathfrak{T}) (1-2s_r)} \label{relevantterms} \end{equation}
as given in Theorem~\ref{thm:bbfm}. However, $\tmmathbf{k}^\Gamma(\mathfrak{T}) = (k_1^\Gamma, \ldots, k_r^\Gamma)$ may be easily recovered from our fixed choice of highest weight $\lambda+\rho=(\ell_1, \cdots, \ell_r, 0)$ and the row sums $d_i := d_i(\mathfrak{T})$ used in the monomial above. Indeed,
$$ k_1^\Gamma = d_1 - (\ell_2 + \cdots + \ell_r) , \quad k_2^\Gamma = d_2 - (\ell_3 + \cdots + \ell_r), \ldots, \quad k_r^\Gamma = d_r. $$
Hence, upon performing this simple transformation, we may recover the monomials in $q^{1-2s_i}$ in~(\ref{relevantterms}) from those in $z_j$ appearing in $\mathcal{G}^{\Gamma} (\mathcal{S},\tmmathbf{z})$ of the above proposition. A similar set of transformations holds for the Delta rules.  
\end{remark}

\begin{proof}
  Proposition~{\ref{matching}} is a consequence of the bijection given in Proposition
  \ref{bijprop}. We sketch the proof for Gamma ice, as the proof for Delta ice
  is similar. Recall that \minussign vertical spins correspond to the entries
  of the pattern, so the values $\gamma (a_{i, j})$ given in (\ref{gammagtwt})
  should appear in the Boltzmann weights for vertices sitting above a \minussign
  vertical spin. The particular cases of (\ref{gammagtwt}) to be used are determined by the
  vertical spin \tmtextit{above} the vertex in question. We now show that each
  $b_{i, j}$ in (\ref{gammagtwt}) matches the charge, the number of \plussign signs
  to the right of the vertex in row $i$. Equivalently, we must show that every
  spin between column $a_{i, j}$ and column $a_{i - 1, j}$ in row $i$ is
  assigned a \plussign. So suppose that $a_{i, j} > a_{i - 1, j}$ and let $v$ be the
  vertex in row $i$, column $a_{i, j}$, and let $w$ be the vertex in row $i$,
  column $a_{i - 1, j}$. Then the north and south spins for $v$ are $(+, -)$
  which, by the six admissible configurations in Gamma ice, forces the east
  spin to be \plussign. All the vertices between $v$ and $w$ have north and south
  spins $(+, +)$ according to our bijection. The east spin \plussign for $v$ becomes
  the west spin for the neighboring vertex $v'$ to the right of $v$, forcing
  the east spin of $v'$ to be \plussign as well. This effect propogates down the
  row, forcing all row spins between $v$ and $w$ to be \plussign. Finally we must
  show that the spectral parameters for $\mathcal{G}$ are given by differences
  of consecutive row sums. This is Lemma 3 of {\cite{hkice}}. \ \ \ \ 
\end{proof}

Given a fixed set of boundary conditions for the vertex model and an
assignment $B$ of Boltzmann weights associated to each admissible
vertex, we refer to the set of all admissible states
$\mathcal{S}$ as a ``system.'' Given a system $\mathfrak{S}$, its
partition function $Z(\mathfrak{S})$ is defined as
\begin{equation} Z(\mathfrak{S}) := Z(\mathfrak{S}, \tmmathbf{z}) =
\sum_{\mathcal{S} \in \mathfrak{S}} B(\mathcal{S}, \tmmathbf{z}),
\quad \text{with} \quad B(\mathcal{S}, \tmmathbf{z}) := \prod_{v \in
\mathcal{S}} B(v), \label{zpartdefined} \end{equation} 
where this latter product is taken over all
vertices $v$ in the state $\mathcal{S}$. In particular, let
$\mathfrak{S}^\Gamma$ denote the system with boundary
conditions as in Section~\ref{bijectionsection}, Boltzmann weights $B^\Gamma$, and
rows labeled in descending order from top to bottom. Similarly, let
$\mathfrak{S}^\Delta$ denote the system with the same boundary
conditions, but with Boltzmann weights $B^\Delta$ and rows labeled in ascending
order from top to bottom. Using this language, we may now
summarize the results of the past two sections in a single theorem.

\begin{theorem} \label{theycommute} Given a dominant weight $\lambda$
for $\tmop{GL}_{r+1}$, the metaplectic Whittaker function
$W(\mathbf{t}_\lambda)$ is expressible as either of the two partition functions
$Z(\mathfrak{S}^\Gamma)$ or $Z(\mathfrak{S}^\Delta).$
\end{theorem}

This is merely the combination of Theorems~\ref{thm:bbfm} and~\ref{statementA} together with Proposition~\ref{matching}.

\section{Transfer Matrices}

Baxter considered the problem of computing partition functions for solvable
lattice models (cf. \cite{baxter}). His approach is based on the idea of using
the Yang-Baxter equation (called the ``star-triangle identity'' by Baxter) to
prove the commutativity of {\tmem{row transfer matrices}}. We will show that
basic properties of metaplectic Whittaker functions can be interpreted as
commutativity of such transfer matrices, and at least when the metaplectic
degree $n = 1$, the Yang-Baxter equation can be used to give proofs of these.

The row transfer matrices shall now be described. Let us consider a row of
vertices that all have the same Boltzmann weights. If
$B = (a_1, a_2, b_1, b_2, c_1, c_2)$ then we use the 
assignment of Boltzmann weights in the following table.
\[ \begin{array}{|c|c|c|c|c|c|c|}
     \hline
     & \includegraphics{gamma1c.mps} & \includegraphics{gamma6c.mps} &
     \includegraphics{gamma4c.mps} & \includegraphics{gamma5c.mps} &
     \includegraphics{gamma2c.mps} & \includegraphics{gamma3c.mps}\\
     \hline
     \begin{array}{c}
       \tmop{Boltzmann}\\
       \tmop{weight}
     \end{array} & a_1 & a_2 & b_1 & b_2 & c_1 & c_2\\
     \hline
   \end{array} \]

The vertical edge spins in the top and bottom boundaries will be collected
into vectors $\alpha = (\alpha_N, \cdots, \alpha_0)$ and $\beta = (\beta_N,
\cdots, \beta_0)$. The subscripts correspond to the columns which, we recall,
are numbered in ascending order from right to left. For example, if $\alpha = (-, +, -, +, +, -)$
and $\beta = (+, -, +, +, +, -)$, we would consider the partition function of
the following one-layer system of ice:
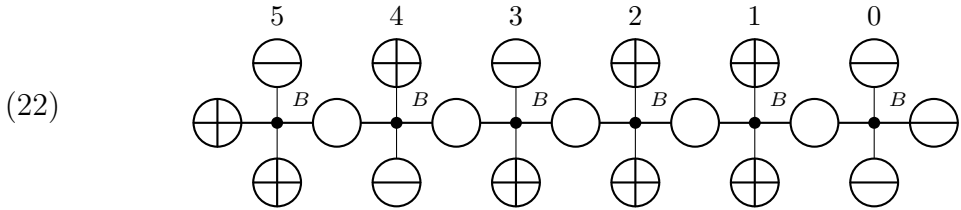
\begin{equation} \vcenter{\hbox to 4 in{\includegraphics[scale=0.9]{transfer.mps}}} 
\label{eq:onelayer} \end{equation}

Let $V_B(\alpha, \beta)$ denote the partition function. Recall that we
compute this as follows. We complete the state by assigning values to
the interior edges (unlabeled in this figure) and sum over all such
completions. Let $V_B$ be the $2^{N + 1} \times 2^{N +
1}$ matrix whose entries are all possible partition functions $V_B(\alpha, \beta)$, 
where the choices of $\alpha$ and $\beta$ index the rows and columns of the matrix,
respectively. This is referred to as the {\it transfer matrix} for the
one-layer system of size $N$ with Boltzmann weights $B$ at every vertex.

Now let us consider a two-layer system:
\begin{equation} \vcenter{\hbox to 4in{\includegraphics[scale=0.9]{transfer1.mps}}} 
\label{eq:twolayer} \end{equation}
Note that we are using two sets of Boltzmann weights $B_1$ and
$B_2$ for the top and bottom layers, respectively.  We may try to
express the partition function $V(\alpha,\gamma)$ for the two layer system
pictured above having top row $\alpha$ as in (\ref{eq:onelayer}) and bottom
row $\gamma = (+, +, +, -, +, +)$ in terms of one-layer partition
functions. However each one-layer system is only determined upon a choice of
vertical spins lying between the rows. One such choice of edge spins is
$\beta$ as in the one-layer system in (\ref{eq:onelayer}), but we must sum
over all possible choices to get the partition function of the two-layer
system. Therefore
\[V(\alpha,\gamma)=\sum_{\beta}
V_{B_1} (\alpha, \beta) V_{B_2} (\beta, \gamma),
\]
which is precisely the entry $V(\alpha, \gamma)$ in the product of the two
transfer matrices $V_{B_1}$ and $V_{B_2}$.

Cases where the transfer matrices commute are of
special interest. Indeed, this commutativity means
that one can interchange the roles of Boltzmann weights 
$B_1$ and $B_2$ in (\ref{eq:twolayer})
and the value of
the product of the transfer matrices is unchanged. Baxter considers
the case where $B_1 = (a, a, b, b, c, c)$ and $B_2 = (a',
a', b', b', c', c')$ for arbitrary choices of $a,a',b,b',c,c'$. However, his boundary 
conditions are toroidal; that is, the
boundary edges at the left and right edges of the each row are identified
and treated as interior edges, hence summed over. With this modification,
Baxter proves that if $\triangle = \triangle'$, where $\triangle = (a^2 + b^2
- c^2) / 2 a b$ and $\triangle'$ is similarly defined with $a', b'$ and $c'$, then the
transfer matrices commute. Obtaining a sufficiently large family of commuting
transfer matrices is a step towards evaluating the partition function, since
by doing so one can make the eigenspaces one-dimensional. Thus the problem of
simultaneously diagonalizing them has a unique solution and therefore becomes
tractable.

Let us now show how Statement~B may be formulated in terms of commuting transfer
matrices. We consider a two layer system having a layer of Gamma ice and a layer of Delta ice,
thus:
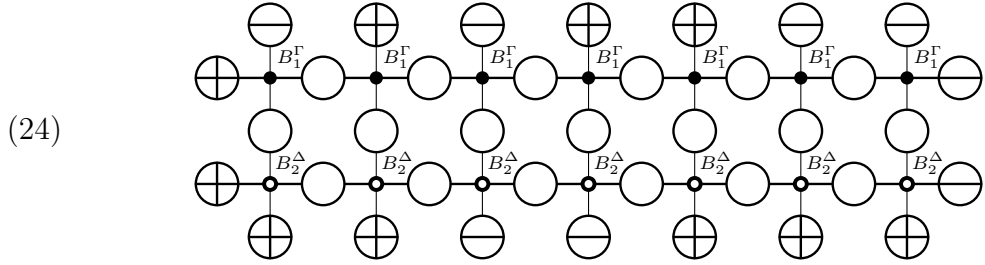
\begin{equation} \vcenter{\hbox to 4in{\includegraphics{transfer2.mps}}} 
\label{gammadeltaice} \end{equation}
We use the values in (\ref{weighttable}); in the top row, the spectral
parameter is $z_1$, and in the bottom row it is $z_2$. Regarding the boundary
conditions, as always the rows of ice must have \plussign at the left edge and
\minussign at the right edge. Furthermore, we fix a choice of spins for
the top edge and the bottom edge of this two-layered ice such that the top edge 
has two more \minussign than the bottom row.
In this example, the locations of the \minussign along the top edge are (reading from
right to left) $0, 1, 4, 6$ and along the bottom edge, they are at $3, 4$. We have
labeled each vertex with $\bullet \, \Gamma_1$ and $\circ \Delta_2$ to remind the
reader of the Boltzmann weights that we are using. We will call this system
$\mathfrak{S}^{\Gamma \Delta}$ and its partition function 
$Z (\mathfrak{S}^{\Gamma \Delta})$.

On the other hand we may consider the same configuration with the roles
of the Boltzmann weights for Gamma and Delta ice switched, as in the figure
below. Note that the boundary conditions remain the same as in (\ref{gammadeltaice}). 
We will refer to this system as $\mathfrak{S}^{\Delta \Gamma}$.
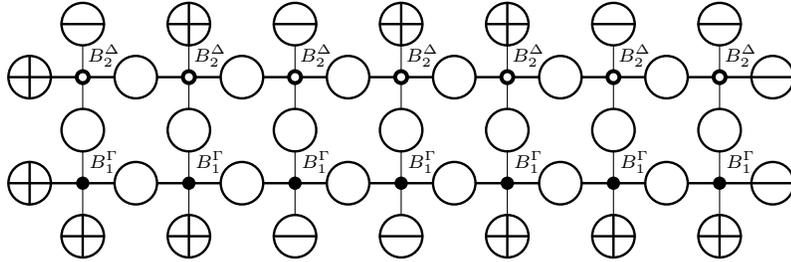
\[ \includegraphics{transfer3.mps} \]

\begin{theorem} \label{abba}
  Given top and bottom boundary values as vectors of spins $\alpha$ and $\gamma$ and Boltzmann
  weights $B_1^{\Gamma}$ and $B_2^{\Delta}$ as in (\ref{weighttable}), let $\mathfrak{S}^{\Gamma \Delta}$ and
  $\mathfrak{S}^{\Delta \Gamma}$ be the systems described above. Then $Z (\mathfrak{S}^{\Gamma \Delta}) 
  = Z(\mathfrak{S}^{\Delta \Gamma})$.
\end{theorem}

We prove this by showing that the claim is equivalent to Statement~B, stated as Theorem~\ref{statementB} here and proved by combinatorial means in \cite{wmd5book}. Note in particular that we have reformulated Statement~B as the commutativity of two transfer matrices.

\medskip

\begin{proof}
  We associate two strictly decreasing vectors of
  integers with $\alpha$ and $\gamma$, which we call $\tmmathbf{l}$ and $\tmmathbf{m}$. Namely, 
  let $\tmmathbf{l} = (l_0, l_1, l_2, \cdots)$,
  where the $l_i$'s are the integers such that $\alpha_{l_i} = -$, arranged in
  descending order; $\tmmathbf{m}$ is defined similarly with regard to $\gamma$. 
  Thus, in the example (\ref{gammadeltaice}) above there are \minussign spins in the $6, 4, 1, 0$ columns of the
  top row and so $\tmmathbf{l} = (6, 4, 1, 0)$, while $\tmmathbf{m} = (4, 3)$. Similarly,
  given any admissible state of the system, let $\beta$ be the middle row of edge
  spins, and associate in similar fashion a sequence $\tmmathbf{a} = (a_1, a_2, \cdots)$
  according to the location of \minussign signs in $\beta$.
  
  We observe that the sequences $\tmmathbf{l}, \tmmathbf{a}, \tmmathbf{m}$ interleave. This holds for the same
  reason that  the rows of the pattern interleave in Proposition~\ref{bijprop};
  it is a consequence of Lemma~2 of {\cite{hkice}}. Therefore the legal
  states of either system $\mathfrak{S}$ are in bijection with the (strict) short Gelfand-Tsetlin
  patterns
  \[ \mathfrak{t =} \left\{ \begin{array}{ccccccccc}
       \ell_0 &  & \ell_1 &  & \cdots &  & \ell_{r - 1} &  & \ell_r\\
       & a_1 &  & a_2 &  & a_{r - 1} &  & a_r & \\
       & & m_1 & & m_2 & & m_{r - 1} & & \end{array} \right\} . \]
These are {\tmem{not}} in bijection with the terms of the sum
$G^{\Gamma \Delta} (\mathfrak{t})$ appearing in
Theorem~\ref{statementB} because there is no condition on the middle
row sum. Rather, the states of ice give all possible middle row sums.
However, letting $\mathcal{G}^{\Gamma \Delta}(\mathcal{S},
\tmmathbf{z})$ denote the Boltzmann weight for a state of $\Gamma
\Delta$ ice, this may be regarded as a homogeneous polynomial in the
two spectral parameters $z_1$ and $z_2$ of our two-row system.  In the
notation of Proposition~\ref{matching}, this monomial is
$z_1^{d_0(\mathfrak{t}) - d_1(\mathfrak{t})} z_2^{d_1(\mathfrak{t}) -
d_2(\mathfrak{t})}$, where $d_i(\mathfrak{t})$ denotes the $i$-th row sum in
the short Gelfand-Tsetlin pattern above. Clearly, the middle row sum can be
recovered from knowledge of this monomial for fixed choice of boundary
conditions $\alpha$ and $\gamma$, which dictate the top and bottom row of the
short pattern.  A similar correspondence may be obtained for the $\Delta
\Gamma$ system whose short patterns $\mathfrak{t}'$ are associated to the
monomial $z_2^{d_0(\mathfrak{t}') - d_1(\mathfrak{t}')}
z_1^{d_1(\mathfrak{t}') - d_2(\mathfrak{t}')}$. Of course, the boundary
conditions remain constant whether we are using the $\Gamma \Delta$ or $\Delta
\Gamma$ system, so $d_0(\mathfrak{t}) = d_0(\mathfrak{t}')$ and
$d_2(\mathfrak{t}) = d_2(\mathfrak{t}')$. Thus, the monomials
$$ z_1^{d_0(\mathfrak{t}) - d_1(\mathfrak{t})} z_2^{d_1(\mathfrak{t}) - d_2(\mathfrak{t})} \quad \text{and} \quad z_2^{d_0(\mathfrak{t}') - d_1(\mathfrak{t}')} z_1^{d_1(\mathfrak{t}') - d_2(\mathfrak{t}')} $$
agree precisely when 
$$ d_1(\mathfrak{t}) = d_0(\mathfrak{t}) + d_2(\mathfrak{t}) - d_1(\mathfrak{t}'), $$
which is exactly the condition on the sum in Theorem~\ref{statementB}. 
Hence we see that the commutativity of transfer matrices -- the statement that
$Z (\mathfrak{S}^{\Gamma \Delta}) = Z(\mathfrak{S}^{\Delta \Gamma})$ -- is an
equality of two homogeneous polynomials and the matching of each monomial
corresponds to the identity of Statement B for each possible middle row sum.
\end{proof}

\section{\label{sec:ybe} The Yang-Baxter Equation} 

The proof of Theorem~\ref{abba}, the commutativity of transfer
matrices, uses the equivalence with Theorem~\ref{statementB} and hence
implicitly relies on all of the combinatorial methods of
\cite{wmd5book} in order to obtain this result. In this section, we
want to explore the extent to which Baxter's methods for solving
statistical lattice models, most notably the Yang-Baxter equation,
may be used to prove the commutativity of transfer matrices.

In our context of two-dimensional square lattice models, the
Yang-Baxter equation may be viewed as a fundamental identity between
partition functions on two very small pieces of ice -- each having 6
boundary edges to be fixed, 3 internal edges, and 3 vertices each with
an assigned set of Boltzmann weights.

\begin{definition}[Yang-Baxter Equation] Let $R,S,$ and $T$ be three
collections of Boltzmann weights associated to each admissible
vertex. Then for every fixed combination of boundary conditions
$\sigma, \tau, \alpha, \beta, \rho, \theta$, we have the following
equality of partition functions:
\begin{equation} \label{startriangleide}
      Z \left( \vcenter{\hbox to
105pt{\includegraphics[scale=0.9]{ybl.mps}}} \right) = Z \left(
\vcenter{\hbox to 105pt{\includegraphics[scale=0.9]{ybr.mps}}} \right)
\; . 
\end{equation} 
Recall that these partition functions are sums of
Boltzmann weights over all admissible states. Hence, the left-hand
side is a sum over all choices of internal edge labels $\mu, \nu,
\gamma$, while the right-hand side is a sum over internal edge labels
$\phi, \psi, \delta$. Note that the roles of $S$ and $T$ are interchanged
on the two sides of the equality.
\end{definition}

In the diagram above one vertex, labeled $R$, has been rotated by
$45^\circ$ for ease of drawing the systems. It should be understood
in the same way as $S$ and $T$ - it has a Boltzmann weight associated
to a set of admissible adjacent edge labels. However, vertices of this
type have a distinguished role to play in the arguments that follow, so
we use the term {\it $R$-vertex} to refer to any vertex rotated by
$45^\circ$ like $R$ in (\ref{startriangleide}).

Once equipped with the Yang-Baxter equation, the commutativity of
transfer matrices, i.e.~invariance of the partition function under
interchange of rows, may be proved under certain assumptions. We
illustrate the method with a three-layer system of ice $\mathfrak{S}$
having boundary conditions and admissible vertices like those of the
system $\mathfrak{S}^\Gamma$, to give the basic idea. Suppose we
wanted to analyze the effect of swapping the second and third rows in
the following configuration:

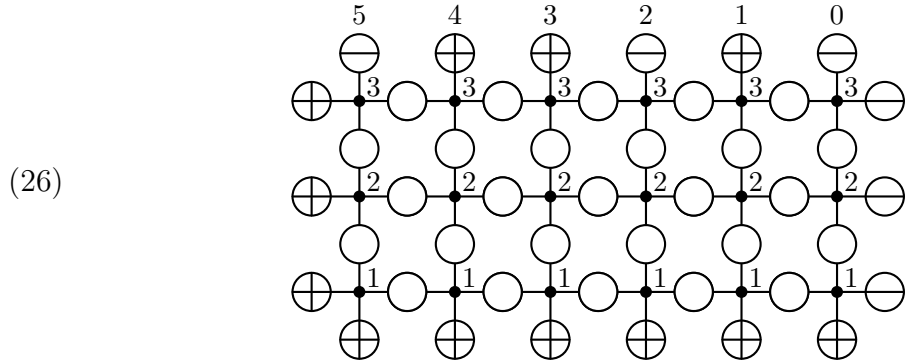
\begin{equation} \label{rectangle}
\vcenter{\hbox to 3in{\includegraphics[scale=0.9]{gamma_ice3.mps}}}
\end{equation}

Suppose there exists only one admissible R-vertex having positive spins on the right;
without loss of generality we take it to have all positive spins. Then
the partition function $Z(\mathfrak{S})$ for~(\ref{rectangle})
multiplied by the Boltzmann weight for the R-vertex with all \plussign
spins is equal to the partition function for the following
configuration of ice.  
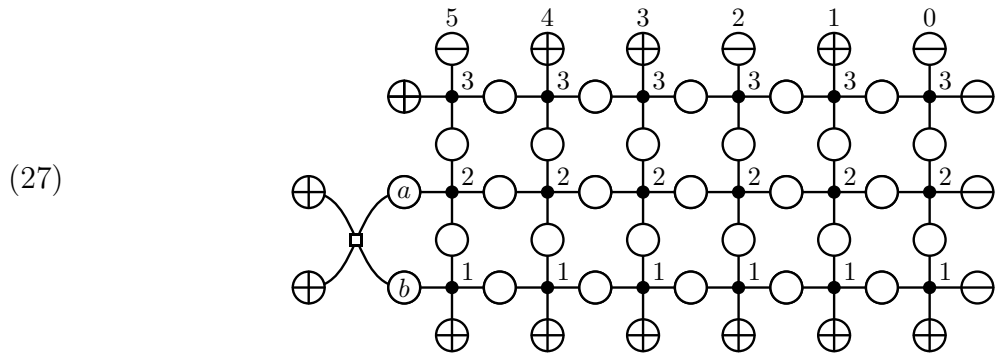
\begin{equation} \vcenter{\hbox to 3in{\includegraphics{state5.mps}}} 
\label{eq:state5} \end{equation}
(By assumption, the only legal
values for $a$ and $b$ are $+$, so every state of this problem
determines a unique state of the original problem.) Now we apply the
Yang-Baxter equation to move this R-vertex rightward, to obtain
equality with the the following configuration. 
 \[ \vcenter{\hbox to 3in
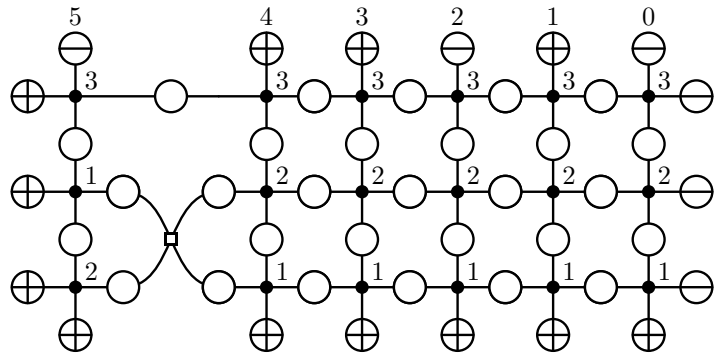
{\includegraphics{state6.mps}}} \] 
 Repeatedly applying the Yang-Baxter equation, we eventually obtain 
 the configuration in which the R-vertex is moved entirely to the right.  
 \begin{equation} \vcenter{\hbox to 3in
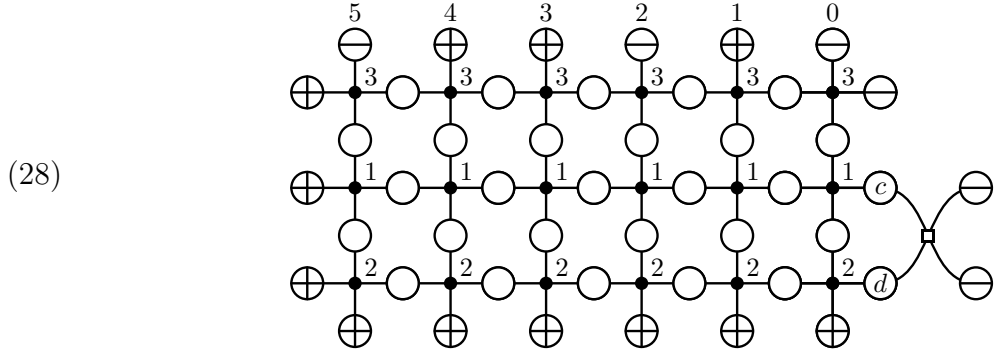
{\includegraphics{state7.mps}}} 
 \label{eq:state7} \end{equation} 
 In drawing the above picture, we have
again assumed that there is just one legal configuration for the
R-vertex having two \minussign spins on the left, and assumed the spins
of this R-vertex were all \minussign. If we let $\mathfrak{S}'$ denote
the system with the same boundary conditions as $\mathfrak{S}$ shown
in (\ref{rectangle}) but with the second and third row Boltzmann
weights interchanged, we have shown
\begin{equation} \label{commute?} B_R
  \left( \vcenter{\hbox to 34pt{\includegraphics{rota1.mps}}} \right)
  Z(\mathfrak{S}) = B_R \left( \vcenter{\hbox to
  34pt{\includegraphics{rota2.mps}}} \right) Z(\mathfrak{S}')
\end{equation}
where $B_R$ denotes the assignment of Boltzmann weight
to each configuration. In particular, if the two admissible
R-vertices coming from the left- and right-hand sides
of~(\ref{commute?}) have equal Boltzmann weights, we obtain the exact
equality of the two configurations, i.e. the commutativity of transfer
matrices.

We now explore the possibility of obtaining a Yang-Baxter equation
with $S$ and $T$ in (\ref{startriangleide}) corresponding to the
Boltzmann weights $B^\Gamma$ and $B^\Delta$, respectively,
from~(\ref{weighttable}).  In light of our previous argument, this
would give an alternate proof of Theorem~\ref{statementB}.  However,
the Boltzmann weights in~(\ref{weighttable}) depend not only on spins
\plussign or \minussign on adjacent edges, but also on a ``charge'' $a$ mod
$n$. Recall from Section~\ref{sec:imwf} that using $B^\Gamma$ weights,
charge records the number of \plussign signs in a row between the given vertex
and the \minussign boundary spin at the right-hand edge the
row. Using $B^\Delta$ weights, charge counts the number of
\minussign signs between the vertex and the \plussign boundary at the left.

In order to demonstrate a Yang-Baxter equation, we need Boltzmann
weights that are purely local --- i.e., depend only on properties of
adjacent edges --- so we need a different way of interpreting
charge. We do this by labeling horizontal edges with both a spin and a
number mod $n$. We declare the Boltzmann weight of these vertices to
be 0 unless the edge labels $a,b$ mod $n$ to the immediate left and right
of the vertex reflect the way charge is counted for the given spins. 
For example, with $B^\Gamma$ weights, $a = b+1$ if the spin below 
$a$ is \plussign and $a=b$ if the spin below $a$
is \minussign. Using this interpretation, we record the non-zero vertices
for both sets of Boltzmann weights:

\begin{equation}
  \begin{array}{|c|l|l|l|l|l|l|l|}
    \hline
    \begin{array}{c}
      \tmop{Gamma}\\
      \text{Ice}
    \end{array} & 
    \vcenter{\vspace{3pt}\hbox{\includegraphics[scale=0.9]{gamma1b.mps}}\vspace{3pt}} &
    \vcenter{\vspace{3pt}\hbox{\includegraphics[scale=0.9]{gamma6b.mps}}\vspace{3pt}} & 
    \vcenter{\vspace{3pt}\hbox{\includegraphics[scale=0.9]{gamma4b.mps}}\vspace{3pt}} &
    \vcenter{\vspace{3pt}\hbox{\includegraphics[scale=0.9]{gamma5b.mps}}\vspace{3pt}} & 
    \vcenter{\vspace{3pt}\hbox{\includegraphics[scale=0.9]{gamma2b.mps}}\vspace{3pt}} &
    \vcenter{\vspace{3pt}\hbox{\includegraphics[scale=0.9]{gamma3b.mps}}\vspace{3pt}}\\
    \hline
    \begin{array}{c}
     B^{\Gamma}  \\
      \text{weight}
    \end{array} 
     & 1 & z_i & g (a) & z_i & h (a) z_i & 1\\
    \hline
    \begin{array}{c}
      \tmop{Delta}\\
      \text{Ice}
    \end{array} & 
    \vcenter{\vspace{3pt}\hbox{\includegraphics[scale=0.9]{delta1b.mps}}\vspace{3pt}} &
    \vcenter{\vspace{3pt}\hbox{\includegraphics[scale=0.9]{delta6b.mps}}\vspace{3pt}} & 
    \vcenter{\vspace{3pt}\hbox{\includegraphics[scale=0.9]{delta4b.mps}}\vspace{3pt}} &
    \vcenter{\vspace{3pt}\hbox{\includegraphics[scale=0.9]{delta5b.mps}}\vspace{3pt}} &
    \vcenter{\vspace{3pt}\hbox{\includegraphics[scale=0.9]{delta2b.mps}}\vspace{3pt}} &
    \vcenter{\vspace{3pt}\hbox{\includegraphics[scale=0.9]{delta3b.mps}}\vspace{3pt}}\\
    \hline
     \begin{array}{c}
      B^{\Delta} \\
      \text{weight} \end{array} & 1 & g(a) z_i & 1 & z_i & 1 & h(a) z_i \\
    \hline
  \end{array}
\label{localweighttable}
\end{equation}
The above vertices are admissible for {\it any} choice of $a$ mod $n$
(and the integers $a+1$ are, of course, understood to be mod $n$ as
well). This means that we are generalizing the six-vertex model, since
due to the dependence on $a$, each vertex has more than six admissible
states.

For $n=1$, the charge labels on horizontal edges are trivial as the
Gauss sums $g(a)$ and $h(a)$ are independent of $a$ as evaluated in
(\ref{nequalone}). For this special case, it was shown in \cite{hkice}
that a Yang-Baxter equation exists with weights $S$ and $T$ as in
(\ref{startriangleide}) taken to be $B^\Gamma$ and $B^\Delta$ from the
table above. We refer the reader to \cite{hkice} for the corresponding
Boltzmann weights $R$ for which the Yang-Baxter equation is
satisfied. Thus we obtain an alternate proof of
Theorem~\ref{theycommute}, or equivalently Theorem~\ref{statementB},
using methods from lattice models.

In general, we know from \cite{wmd5book} that
Theorem~\ref{theycommute} is true for any positive integer $n$.  It
would be extremely interesting to find a local relation like
(\ref{startriangleide}) similarly proving that the transfer matrices
commute, and this is currently under investigation by the authors.

\section{Weyl group invariance and the Yang-Baxter equation}

\newcommand{\SSS}{\mathcal{S}}
\newcommand{\zzz}{\tmmathbf{z}}
\newcommand{\fS}{\mathfrak{S}}

Kazhdan and Patterson \cite[Lemma 1.3.3]{kp} describe how the metaplectic Whittaker
functions transform under the action of the Weyl group.  This
invariance---which does not follow directly from the description
of the coefficients $H$ given in Theorem~\ref{thm:bbfm}---plays a
key role in the proof of the metaplectic Casselman-Shalika formula for
$GL_{r+1}$ by Chinta and Offen \cite{ChintaOffen}, and was the main
inspiration for the Weyl group action in \cite{CG}.

In this section we restate this Weyl group invariance in terms of the
partition functions defined in the previous sections.  We content
ourselves to describe how a simple reflection acts on the partition
function.  Let $\sigma_i$ denote the simple reflection in the Weyl
group corresponding to the $i$-th simple root.  We let $\sigma_i$ act
on the spectral parameter $\tmmathbf {z}=(z_1,z_2, \ldots ,z_{r+1})$
by $\sigma_i(\zzz)=(z_1, \ldots, z_{i-1}, z_{i+1}, z_i,
z_{i+2},\ldots, z_{r+1}),$ i.e. the $i^{th}$ and $(i+1)^{st}$
coordinates are transposed.  Here the notation $Z(\mathfrak S, \zzz)$
refers to the partition function associated to the system
$\mathfrak S$, where $\mathfrak S$ is either of the two systems
$\mathfrak S^\Gamma$ or $\mathfrak S^\Delta$ introduced in
Section~\ref{sec:imwf}.

Further define, for $j=0,\dotsc ,n-1$,
\begin{equation}
   \label{defn:PQ}
   P^{(j)}(x,y)=x^jy^{n-j}\frac {1-q^{-1}}{x^n-q^{-1}y^n}
\mbox{\ \ and\ \ \ } Q^{(j)}(x,y)=g(j)
\frac{x^n-y^n}{x^n-q^{-1}y^n}, 
\end{equation}
where we again use the shorthand notation of (\ref{gandhforshort})
and interpret $g(0) := g(n) =-q^{-1}.$  The functions $P$,$Q$ are
closely related to the functions $\tau^{1}_{s}$,$\tau^{2}_{s}$ of
\cite[Lemma 1.3.3]{kp}.

For each $1\leq i \leq r$, we may decompose the partition function
\begin{equation}\label{eq:decomposition}
Z(\mathfrak{S}, \zzz) =\sum_{0\leq j <n}Z_i^{(j)}(\fS,\zzz),
\end{equation}
where $Z_i^{(j)}(\fS,\zzz)$ is the sum over all states $\SSS \in \fS$ such that
$B(\SSS, \zzz)$ is equal to a constant times
$z_1^{a_1}\cdots z_{r+1}^{a_{r+1}}$ 
where $a_{i}-a_{i+1} \equiv j \; (\text{mod } n)$.  
Then the Whittaker function satisfies
\begin{equation}
   \label{eq:fe} Z_i^{(j)} (\fS , \sigma_{i}(\zzz) ) = P^{(j)}(z_{i+1}, z_{i})\cdot
Z_i^{(j)}(\fS , \zzz) + Q^{(j)}(z_{i+1}, z_{i})\cdot
Z_i^{(n-j)}(\fS , \zzz).
\end{equation}

We now consider the extent to which the functional equations
\eqref{eq:fe} can be interpreted in the language of transfer matrices.  
First we consider the case $n=1$.  
The decomposition on the right-hand side of \eqref{eq:decomposition} 
has only one term, namely $Z$ itself, since
the congruence condition is automatically satisfied by all monomials
for any $i$.  The $i$-th functional equation \eqref{eq:fe} becomes
\[
Z (\fS , \sigma_{i}(\zzz)) = (P^{(0)} (z_{i+1},z_{i})+Q^{(0)} (z_{i+1},z_{i})) Z
(\fS , \zzz ),
\]
or 
\begin{equation}\label{eq:commutation}
(z_{i}-z_{i+1}/q)Z (\fS , \zzz) = (-z_{i}/q + z_{i+1})Z (\fS ,\sigma_{i} (\zzz)).
\end{equation}
Recalling the effect of $\sigma_i$ on $\zzz$ defined above, the
partition function on the right-hand side is the result of swapping
the spectral parameters associated to rows $i$ and $i+1$ in the system
$\fS$. Note that \eqref{eq:commutation} is not exactly the same as
``commutation of two transfer matrices'' because we do not have the
identity $Z (\fS , \zzz) = Z (\fS , \sigma_{i} (\zzz))$.  Indeed, the
partition function $Z$ is not a symmetric function, but it is very
close to one: it is a Schur polynomial times a $q$-deformation of the
Weyl denominator (cf.~\cite{hkice}).

Nevertheless, with assumptions as in Section~\ref{sec:ybe}, we may ask
for a Yang-Baxter equation leading to a proof of
\eqref{eq:commutation}. That is, we seek sets of Boltzmann weights
$R,S,$ and $T$ satisfying~\eqref{startriangleide} where $S = T =
B^\Gamma$ or $S = T = B^\Delta$.  Comparing \eqref{eq:commutation}
with \eqref{commute?}, we further require Boltzmann weights $B_R$ for
the R-vertices such that
\begin{align*} B_R
\left( \vcenter{\hbox to 34pt{\includegraphics{rota1.mps}}} \right)
 &= z_{i}-z_{i+1}/q, \\
B_R \left( \vcenter{\hbox to
34pt{\includegraphics{rota2.mps}}} \right) &= -z_{i}/q + z_{i+1}.
\end{align*}
It follows from results in~\cite{hkice} that we may use the following
coefficients in the R-vertex for Gamma ice:
\begin{equation}
\label{gammagammaweights}
\begin{array}{|c|c|c|c|c|c|}
     \hline
     \includegraphics{rota1.mps} & \includegraphics{rota2.mps} &
     \includegraphics{rotb1.mps} & \includegraphics{rotb2.mps} &
     \includegraphics{rotc1.mps} & \includegraphics{rotc2.mps}\\
     \hline
     z_i - q^{- 1} z_{i + 1} & z_{i + 1} - q^{- 1} z_i & q^{- 1} (z_{i + 1} -
     z_i) & z_{i + 1} - z_i & (1 - q^{- 1}) z_{i + 1} & (1 - q^{- 1}) z_i\\
     \hline
   \end{array}
\end{equation}
We are taking all $t_i=-q^{-1}$ in Table~1 in~\cite{hkice}, and
observe that the order of the rows in this paper are opposite those
in that paper. Our convention here is the same as in~\cite{wmd5book}.

For $n>1$ the situation is more complicated, but rather suggestive.
In general $Z_{i}^{(j)}(\fS) \not = Z_{i}^{(n-j)}(\fS)$, so we cannot
rewrite \eqref{eq:fe} to look like \eqref{commute?} and
\eqref{eq:commutation}.  However according to \eqref{defn:PQ}, the
denominators of $P$ and $Q$ appearing in the $i$th functional equation
\eqref{eq:fe} are equal and independent of $j$. For any $j$, they are
$z_{i}^{n}-z_{i+1}^{n}/q$.  Thus, clearing denominators, we may rewrite
\eqref{eq:fe} as follows
\begin{eqnarray}
  \label{fewoden}
  (z_{i + 1}^n - q^{- 1} z_i^n) Z_i^{(j)} (\mathfrak{S}, \sigma_i
  (\tmmathbf{z})) & = & \nonumber\\
  p^{(j)} (z_{i + 1}, z_i) \cdot Z^{(j)}_i (\mathfrak{S}, \tmmathbf{z}) +
  q^{(j)} (z_{i + 1}, z_i) \cdot Z^{(n - j)}_i (\mathfrak{S}, \tmmathbf{z}) & 
  & 
\end{eqnarray}
where
\[ p^{(j)} (z_{i + 1}, z_i) = (1 - q^{- 1}) z_{i + 1}^j z_i^{n - j},
   \hspace{2em} q^{(j)} (z_{i + 1}, z_i) = g (j) (z_{i + 1}^n - z_i^n) . \]

Let $\mathcal{S}$ be a state of the system, and as before let
$a_1,\cdots,a_{r+1}$ be the exponents of $z_1,\cdots,z_{r+1}$ 
in $B (\mathcal{S}, \tmmathbf{z})$. We make the following observation. In the
weights (\ref{localweighttable}), there is a contribution of $z_i$ if and only
if the charge is not augmented as we move across the vertex. Since (in Gamma
ice) the charges at the right edge will have value 0, it follows that the
charges at the left edge will have value $c_i$ where $a_i + c_i$ is the number
of vertices in the row. Therefore
\begin{equation}
  \label{constantdiff} a_i - a_{i + 1} = c_{i + 1} - c_i
\end{equation}
and we may therefore write
\[ Z_i^{(j)} (\mathfrak{S}, \tmmathbf{z}) = \sum_{\text{$c_{i + 1} - c_i \equiv j$
   mod $n$}} B (\mathcal{S}, \tmmathbf{z}) . \]

We will now explain how, with a suitable R-vertex, (\ref{fewoden}) could also
be interpreted as an identity similar to (\ref{commute?}), but now with sets
of Boltzmann weights involving charges. We will describe the characteristics
that such an R-vertex might have. For simplicity, we will assume that $n$ is
odd.

The value will depend on the spins and charges of the adjacent edges. Let us
assume first that the spins on these four edges are all $+$, with charges
$d_{i + 1}$, $d_i$, $d_{i + 1}'$, $d_i'$ as follows:
\begin{equation}
  \label{pprvertex} \includegraphics{rota1a.mps}
\end{equation}
If $j = d_{i + 1} - d_i$ and $j' = d_{i + 1}' - d_i'$ then we require that the
Boltzmann weight of this vertex $v$ is zero unless $j' \equiv j$ or $n - j$
mod $n$. Moreover in these cases we require that the Boltzmann weight of
(\ref{pprvertex}) is
\[ \left\{ \begin{array}{ll}
     p^{(j)} (z_{i + 1}, z_i) & \text{if $j \equiv j'$ mod $n$} \\
     q^{(j)} (z_{i + 1}, z_i) & \text{if $j \equiv n - j'$ mod $n$}
     \end{array}
\right.\]
except when $j \equiv 0$. In this case the weight will be
\[ p^{(0)} (z_{i + 1}, z_i) + q^{(0)} (z_{i + 1}, z_i) = z_i^n - q^{- 1} z_{i
   + 1}^n, \]
since $g (0) = - q^{- 1}$.

Regarding the case where the vertex has spin $-$ on all four adjoining edges,
we require that the Boltzmann weight of
\[ \includegraphics{rota2a.mps} \]
is zero unless $d_i = d_{i + 1} = 0$, in which case it is $z_{i + 1}^n - q^{-
1} z_i^n$.

Assuming that the R-vertex has the above properties, we may now express the
functional equation in a form similar to (\ref{commute?}). Let us fix the
vertical edge spins above the $z_{i + 1}$ row and below the $z_i$ row, and
work with just the two relevant rows; let $\mathfrak{S}'$ denote the
two-layer system consisting of just rows $i+1$ and $i$ with these boundary
spins fixed. In order to establish (\ref{eq:fe}), or equivalently
(\ref{fewoden}), it suffices to show
\begin{eqnarray*}
  (z_{i + 1}^n - q^{- 1} z_i^n) Z_i^{(j)} (\mathfrak{S}', \sigma_i
  (\tmmathbf{z})) & = & \\
  p^{(j)} (z_{i + 1}, z_i) \cdot Z^{(j)}_i (\mathfrak{S}', \tmmathbf{z}) +
  q^{(j)} (z_{i + 1}, z_i) \cdot Z^{(n - j)}_i (\mathfrak{S}', \tmmathbf{z}). & 
  & 
\end{eqnarray*}
Since $Z(\mathfrak{S}',\mathbf{z})$ is a homogeneous polynomial in the $z_i$,
and since only $a_i$ and $a_{i + 1}$ are allowed to vary, we have $a_i + a_{i+ 1}$ 
equal to a constant. Since we are assuming that $n$ is odd, there will be a
unique pair of charges $c_i$ and $c_{i + 1}$ mod $n$ such that
(\ref{constantdiff}) is satisfied, and such that $c_{i+1}-c_i\equiv j$
modulo~$n$.

Now let us consider the partition function of the system
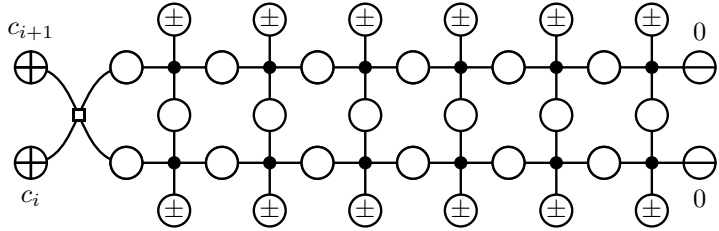
\[ \includegraphics{state8.mps} \]
obtained by attaching the R-vertex to the left of $\mathfrak{S}'$.  From the
above discussion, this equals
\[ p^{(j)} (z_{i + 1}, z_i) \cdot Z^{(j)}_i (\mathfrak{S}', \tmmathbf{z}) +
   q^{(j)} (z_{i + 1}, z_i) \cdot Z^{(n - j)}_i (\mathfrak{S}', \tmmathbf{z}).
\]
Similarly the partition function of the system
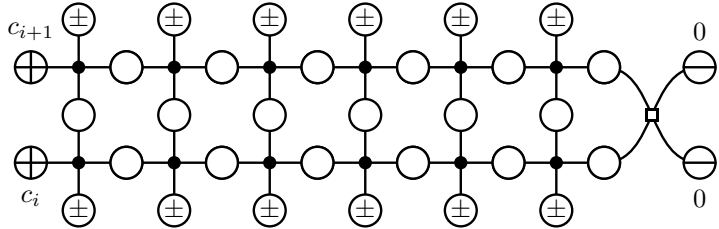
\[ \includegraphics{state9.mps} \]
is
\[ (z_{i + 1}^n - q^{- 1} z_i^n) Z_i^{(j)} (\mathfrak{S}', \sigma_i
   (\tmmathbf{z})) . \]
The equality of these partition functions implies (\ref{fewoden}).

At this writing, we do not know if the values of the R-vertex that we have
described can be completed to a full R-matrix such that the appropriate
Yang-Baxter equation is satisfied. We know that this can be done when $n = 1$,
and since (\ref{fewoden}) is true, it seems very plausible that this
can be done in general. Thus we may conjecture that within this scheme, or
some similar one, it is possible to formulate a Yang-Baxter equation adapted
to these weights that gives a proof of (\ref{fewoden}). Such a ``metaplectic''
Yang-Baxter equation might well have importance beyond the problems that we
have discussed in this paper.

\bibliographystyle{amsplain_initials}
\bibliography{mice}

\def\cprime{$'$}
\providecommand{\bysame}{\leavevmode\hbox to3em{\hrulefill}\thinspace}
\providecommand{\MR}{\relax\ifhmode\unskip\space\fi MR }
\providecommand{\MRhref}[2]{%
  \href{http://www.ams.org/mathscinet-getitem?mr=#1}{#2}
}
\providecommand{\href}[2]{#2}
\begin{thebibliography}{10}

\bibitem{BanksLevySepanski}
W.~D. Banks, J.~Levy, and M.~R. Sepanski, \emph{Block-compatible metaplectic
  cocycles}, J. Reine Angew. Math. \textbf{507} (1999), 131--163.

\bibitem{baxter}
R.~J. Baxter, \emph{Exactly solved models in statistical mechanics}, Academic
  Press Inc. [Harcourt Brace Jovanovich Publishers], London, 1989, Reprint of
  the 1982 original.

\bibitem{BerensteinZelevinsky}
A.~Berenstein and A.~Zelevinsky, \emph{Tensor product multiplicities, canonical
  bases and totally positive varieties}, Invent. Math. \textbf{143} (2001),
  no.~1, 77--128.

\bibitem{hkice}
B.~Brubaker, D.~Bump, and S.~Friedberg, \emph{{Schur polynomials and the
  Yang-Baxter equation}}, submitted.

\bibitem{eisenxtal}
\bysame, \emph{{Weyl group multiple {D}irichlet series, {E}isenstein series and
  crystal bases}}, to appear in Annals of Math.

\bibitem{wmd5book}
\bysame, \emph{{Weyl Group Multiple {D}irichlet Series: {T}ype {A}
  Combinatorial Theory}}, to appear in Annals of Mathematics Studies.

\bibitem{wmd3}
B.~Brubaker, D.~Bump, S.~Friedberg, and J.~Hoffstein, \emph{Weyl group multiple
  {D}irichlet series. {III}. {E}isenstein series and twisted unstable {$A\sb
  r$}}, Ann. of Math. (2) \textbf{166} (2007), no.~1, 293--316.

\bibitem{gelbvol}
B.~Brubaker, D.~Bump, and S.~Friedberg, \emph{Gauss sum combinatorics and
  metaplectic {E}isenstein series}, Automorphic forms and {$L$}-functions {I}.
  {G}lobal aspects, Contemp. Math., vol. 488, Amer. Math. Soc., Providence, RI,
  2009, pp.~61--81.

\bibitem{ChintaOffen}
G.~Chinta and O.~Offen, \emph{A metaplectic {Casselmann--Shalika} formula for
  {$GL_r$}}, submitted.

\bibitem{CG}
G.~Chinta and P.~E. Gunnells, \emph{Constructing {W}eyl group multiple
  {D}irichlet series}, J. Amer. Math. Soc. \textbf{23} (2010), no.~1, 189--215.

\bibitem{Friedberg-McNamara}
S.~Friedberg and P.~McNamara, in preparation.

\bibitem{hamelkinguturn}
A.~M. Hamel and R.~C. King, \emph{U-turn alternating sign matrices, symplectic
  shifted tableaux and their weighted enumeration}, J. Algebraic Combin.
  \textbf{21} (2005), no.~4, 395--421.

\bibitem{kp}
D.~A. Kazhdan and S.~J. Patterson, \emph{Metaplectic forms}, Inst. Hautes
  \'Etudes Sci. Publ. Math. (1984), no.~59, 35--142.

\bibitem{BerensteinKirillov}
A.~N. Kirillov and A.~D. Berenstein, \emph{Groups generated by involutions,
  {G}el\cprime fand-{T}setlin patterns, and combinatorics of {Y}oung tableaux},
  Algebra i Analiz \textbf{7} (1995), no.~1, 92--152.

\bibitem{Littelmann}
P.~Littelmann, \emph{Cones, crystals, and patterns}, Transform. Groups
  \textbf{3} (1998), no.~2, 145--179.

\bibitem{Matsumoto}
H.~Matsumoto, \emph{Sur les sous-groupes arithm\'etiques des groupes
  semi-simples d\'eploy\'es}, Ann. Sci. \'Ecole Norm. Sup. (4) \textbf{2}
  (1969), 1--62.

\bibitem{McNamara}
P.~McNamara, \emph{Metaplectic {W}hittaker functions and crystal bases}, to
  appear in Duke Math. J.

\bibitem{mcnamarafoundations}
\bysame, \emph{Principal series representations of metaplectic groups over
  local fields}, this volume.

\bibitem{Shintani}
T.~Shintani, \emph{On an explicit formula for class-{$1$} ``{W}hittaker
  functions'' on {$GL_{n}$} over {$P$}-adic fields}, Proc. Japan Acad.
  \textbf{52} (1976), no.~4, 180--182.

\end{thebibliography}

\end{document}